\documentclass[3p,preprint]{elsarticle}

\usepackage{amsmath}
\usepackage{amsfonts}
\usepackage{amssymb,latexsym}
\usepackage{amsthm}
\usepackage{graphicx}

\usepackage{color}
\usepackage{newlfont} 
\usepackage{subfigure}
\usepackage{verbatim}
\usepackage{rotating}
\usepackage{multirow}
\usepackage{units}
\usepackage{bm}
\usepackage[english]{babel}
\usepackage[utf8]{inputenc}
\usepackage{algorithm}
\usepackage[noend]{algpseudocode} 
\usepackage{booktabs}
\usepackage{caption}
\usepackage{hyperref}
\usepackage{subfigure}
\usepackage[export]{adjustbox}

\journal{...}

\begin{document}

\begin{frontmatter}

\title{Comparing algebraic cubature rules\\ on spline curved elements}



\author[address-PD]{Alvise Sommariva}
\ead{alvise@math.unipd.it}

\author[address-PD]{Marco Vianello}
\ead{marcov@math.unipd.it}

\cortext[corrauthor]{Corresponding author: Marco Vianello}

\address[address-PD]{Dept. of Mathematics, University of Padova, via Trieste 63, 35121 Padova (Italy)}

\begin{abstract}
We compare four methods for the construction of algebraic cubature rules on planar elements, whose boundary is tracked by splines. The methods, that we have developed over the last two decades, are based on  Green theorem together with some cornerstones of polynomial approximation theory: Gaussian quadrature, Tchakaloff theorem, discretized Chebyshev expansion (hyperinterpolation), Fekete-like interpolation. We discuss their advantages and drawbacks in view of the application to curved polytopal element methods. We have also made freely available at a single site the corresponding open-source Matlab codes. 
\end{abstract}


\end{frontmatter}

\section{Introduction}
The availability of efficient cubature methods is a relevant issue within FEM and VEM methods on nonstandard  elements, especially for curved polytopal elements in high-order instances. Indeed, computation of stiffness and mass matrices can still represent a bottleneck, that has to be coped with by specialized cubature rules. During the last two decades, a devoted literature has been growing on this topic, with some attention on avoiding subtriangulations or in general subtessellations; with no pretence of exhaustivity, we may quote \cite{AHP18,ABD20,ASV20,CDMS24,CLS15,CS21,D26,GWE20,KZPGMBdVD,LvPBD24,MSSLP24,MXS10,RSV26,SFMH11,SV09,SV21} and the references therein.
  
In this respect, the main purpose of this paper is to give a guided tour in the use of some specific codes, that we developed over the last twenty years, for algebraic cubature on a general {\em curved element} $E\subset \mathbb{R}^2$ (either convex or concave), in the case when its {\em piecewise regular Jordan boundary $\partial E$ is tracked by splines}. Such codes  produce weights $\{w_s\}$ and nodes $\{P_s=(x_s,y_s)\}$ such that 
\begin{equation} 
\int_E{f(x,y)\,dx dy}=\sum_{s=1}^M{w_s\,f(P_s)}\;\;\forall f\in \mathbb{P}_n\;,
\end{equation}
(where $\mathbb{P}_n$ denotes the space of bivariate polynomials with total degree not exceeding $n$), and are here named correspondingly to the construction approach:

\begin{itemize}
\item {\bf Spline-Gauss} (2007-2009) \cite{SV07,SV09}: collection of product-like Gaussian rules by Green theorem and numerical primitives ({\em possible negative weights and external nodes on concave elements, high cardinality, empirical stability in general and theoretical stability on convex elements, low construction cost}).

\item {\bf Spline-Tchakaloff} (2021 \& 2024) \cite{SV21,SV24}: moment-matching by Green theorem and NonNegative Least Squares, supported at Tchakaloff sets ({\em subtessellation-free, positive weights and internal nodes, low cardinality, theoretical stability, high construction cost}).

\item {\bf Spline-Chebyshev} (2025-2026) \cite{RSV26,RSV26-2,SV25}: cheap moment-matching by Green theorem and discretized  Chebyshev orthogonal expansion (hyperinterpolation), supported at a minimal cubature set for the Chebyshev measure in a bounding box ({\em subtessellation-free, some external nodes and possible negative weights, low cardinality, theoretical stability, low construction cost}).

\item {\bf Spline-Fekete} (2009 \& 2026) \cite{SV09-2}: moment-matching by Green theorem and QR factorization with column pivoting, supported at approximate Fekete sets ({\em subtessellation-free, possible negative weights but internal nodes, low cardinality, empirical stability, moderate  construction cost}).

\end{itemize}
By stability of the cubature rule we mean the fact that $\|\mathbf{w}\|_1$ is bounded with $n$, and more precisely that, also in the presence of some negative weights, the stability parameter 
\begin{equation} \label{sigma}
\sigma=\frac{\sum_s{|w_s|}}{\sum_s{w_s}}=\frac{\sum_s{|w_s|}}{area(E)}\geq 1   
\end{equation} 
stays relatively close to   
the optimal value 1 (positive weights), at least for the degrees of interest 
in FEM and VEM applications. 

In the sequel we give 
a brief description of the underlying method for each code, pointing out its theoretical foundations, as well as  advantages and drawbacks. Moreover, we illustrate the performance on some test elements. The corresponding open-source codes are implemented in Matlab and freely available at \cite{CUB26}.

\section{The four cubature rules}

The starting point of all the methods is that the boundary $\partial E$ of the element is piecewise regular, and each regular arc is parametrized by a number of sampling points with a spline curve of local degree $d_i$. Namely
$$
\partial E=\bigcup_{i=1}^m{\Gamma_i}\;,\;\Gamma_i=\bigcup_{j=1}^{s_i}{\Gamma_{i,j}}\;,
\;\;\Gamma_{i,j}=Q_{i,j-1}\frown Q_{i,j}\;,
$$
$$
\{Q_{i,j}\}\subset \Gamma_i=S_{i}([0,1])\;,\;\;S_{i}(t)=(x_{i}(t),y_{i}(t))\;,\;t\in [0,1]\;,
$$
\begin{equation} \label{splineb}
Q_{i,j}=S_{i}(t_{i,j})\;,\;0=t_{i,0}<t_{i,1}<\dots<t_{i,s_i}=1\;,
\end{equation}
where $S_i=S_{i,d_i}\in \mathbb{P}_{d_i}$ on each parameter subinterval $[t_{i,j-1},t_{i,j}]$ (for simplicity we consider a fixed parametrization interval $[0,1]$ and a fixed spline degree $d$). The joining points $\{Q_{i,s_i}\}$ of two consecutive regular arcs are the ``vertices'' of such a spline curved polygon. In practice, the code user can simply give the (counterclockwise) ordered point subsequences $\{Q_{i,j}\}$ for each $i$ together with the required spline degree $d$, and the spline boundary is automatically set up. 

All the methods use in a substantial way Green theorem (cf. e.g. \cite{A69}) to compute integrals of polynomials, observing that a  bivariate polynomial $f(x,y)$ of degree $n$ restricted to a spline boundary, is a univariate polynomial on each 
polynomial spline arc $\Gamma_{i,j}$, namely the polynomial $f(x_{i}(t),y_{i}(t))$ of degree $dn$ in the $t$ variable. Now, if for example an $x$-primitive 
\begin{equation} \label{prim}
\mathcal{F}(x,y)=\int{f(x,y)\,dx}
\end{equation}
of a polynomial $f(x,y)\in \mathbb{P}_n$ is known then 
$$
\int_E{f(x,y)\;dxdy}=\oint_{\partial E}{\mathcal{F}(x,y)\,dy}=\sum_{i=1}^m\int_{\Gamma_i}{\mathcal{F}(x,y)\,dy}
$$
$$
=\sum_{i=1}^m\sum_{j=1}^{s_i}\int_{\Gamma_{i,j}}{\mathcal{F}(x,y)\,dy}
=\sum_{i=1}^m\sum_{j=1}^{s_i}\int_{t_{i,j-1}}^{t_{i,j}}{\mathcal{F}(S_{i}(t))\,y'_{i}(t)\;dt}
$$

\begin{equation} \label{splinegauss}
=\sum_{i=1}^m\sum_{j=1}^{s_i}\sum_{h=1}^{\lfloor \frac{(n+2)d-1}{2}\rfloor}{\lambda _{i,j,h}\,\mathcal{F}(S_{i}(\tau_{i,j,h}))\,y'_{i}(\tau_{i,j,h})}\;,
\end{equation}
where $\{\tau_{i,j,h}\}_h$ and $\{\lambda _{i,j,h}\}_h$ are the nodes and positive weights of the Gauss-Legendre formula with exactness degree $(n+2)d-1$ on the interval $[t_{i,j-1},t_{i,j}]$. Indeed, notice that the degree of $\mathcal{F}(S_{i}(t))$ is there at most $(n+1)d$ since $\mathcal{F}$ is an $x$-primitive of a polynomial in $\mathbb{P}_n$, and the degree of $y'_{i}(t)$ is at most $d-1$.

\subsection{{\bf Spline-Gauss}}
The original algorithm \cite{SV09}, which stemmed from a previous work on cubature for linear polygons \cite{SV07}, has been 
widely used in the finite and virtual elements literature over the last decade; cf. e.g. \cite{D26,KZPGMBdVD}, only to quote two very recent ones in the VEM framework. Though the method doesn't apparently make a preliminary subtessellation of the element, it is ultimately based on a sort of implicit covering by disjoint curved trapezoids and piecewise product Gauss-Legendre quadrature. 

It works substantially as follows: 
\begin{itemize}
\item a vertical reference line $x=a$ is fixed (in practice, any line up to a rotation)

\item a primitive $\mathcal{F}$ is computed for every $i,j,h$ at the boundary points $S_{i}(\tau_{i,j,h})=(x_{i,j,h},y_{i,j,h})\,$ $=(x_{i}(\tau_{i,j,h}),y_{i}(\tau_{i,j,h}))$
in (\ref{splinegauss}) as 
$$
\mathcal{F}(S_{i}(\tau_{i,j,h}))=\int_a^{x_{i,j,h}}{f(x,y_{i,j,h})\,dx}=\sum_{k=1}^{\lfloor n/2\rfloor}{\nu_{i,j,h,k}\,f(P_{i,j,h,k})}\;,\;\forall f\in \mathbb{P}_n\;,
$$
with $P_{i,j,h,k}=(\xi_{i,j,h,k},y_{i,j,h})$,  where $\{\xi_{i,j,h,k}\}_k$ and $\{\nu_{i,j,h,k}\}_k$ are the nodes and weights of the Gauss-Legendre rule of exactness degree $dn$ on the interval $[a,x_{i,j,h}]$.

\item the final algebraic cubature rule on the element is 
$$
\int_E{f(x,y)\;dxdy}=\sum_{i=1}^m\sum_{j=1}^{s_i}\sum_{h=1}^{\lfloor \frac{d_in+2d_i-1}{2}\rfloor}{\sum_{k=1}^{\lfloor n/2\rfloor}{w_{i,j,h,k}\,f(P_{i,j,h,k})}}\;,\;\forall f\in \mathbb{P}_n\;,
$$
\begin{equation} \label{splineG}
w_{i,j,h,k}=\lambda _{i,j,h}\,y'_{i}(\tau_{i,j,h})\,\nu_{i,j,h,k}\;. 
\end{equation}

\end{itemize}

Though the construction of the Spline-Gauss cubature rule (\ref{splineG}) is easy and cheap, being based on univariate Gauss-Legendre quadrature, its {\em cardinality can be large}, namely
\begin{equation} \label{cardsplineG}
M=card(\{P_{i,j,h,k}\})\approx \frac{n^2}{4}\,\sum_{i=1}^{m}{d_is_i}\;,
\end{equation}
so that $M$  can already reach a size of the thousands for example with cubic splines ($d_i=3$) and relatively small values of $n$ (exactness degree) and $m$ (number of regular sides of the curved element).

In addition, the Spline-Gauss rule has in general {\em possible external nodes} and some {\em possible negative weights}. This happens typically with concave elements or with nonoptimal choices of the reference line. Nevertheless, as proved in \cite{SV09}, the stability parameter $\sigma$ in (\ref{sigma}) is bounded. Moreover in practice, at least for moderate  exactness degrees (up to around $n=20$) the negative weights turn out to be few and of small size, so that the stability parameter stays relatively close to the optimal value 1. On the other hand, again in  \cite{SV09} it is proved that with {\em convex} spline boundaries the nodes result all {\em internal} and the weights all {\em nonnegative}, simply by choosing as reference line a diameter of the element.

It is worth recalling that, when the weights of Spline-Gauss are positive, one could compress the corresponding discrete measure obtaining a Tchakaloff-like rule \cite{T57}, by the NonNegative Least-Squares method proposed 
in \cite{SV15}, and then adopted for example in \cite{D26,KZPGMBdVD} in the VEM framework. However, we prefer to describe in the next subsection a method, 
named Spline-Tchakaloff, that produces a Tchakaloff-like rule with guaranteed internal nodes and positive weights, irrespectively of the element shape.

\subsection{{\bf Spline-Tchakaloff}}
The method relies on {\em measure compression via linear or quadratic programming}, a concept emerged in the numerical literature around the mid 2010s, independently in different deterministic and probabilistic contexts; cf. e.g. \cite{LL12,MP13,RB15,SV15,Tche15} with the references therein.
The theoretical background is the celebrated Tchakaloff theorem \cite{T57}, a cornerstone of cubature theory, which substantially asserts that 
\begin{itemize}
\item 
{\em there exists an exact cubature formula with at most $N=dim(\mathbb{P}_n)$ nodes and positive weights, for the integration of polynomials in $\mathbb{P}_n$ with respect to a positive Borel measure with compact support}. 
\end{itemize}
An often overlooked refinement of Tchakaloff theorem is a general result 
due to Wilhelmsen \cite{W76}, following a previous construction by Davis \cite{D67}, which ensures that Tchakaloff-like rules on a compact domain $E$ can be supported at a suitable subset of a so-called Tchakaloff set (a sufficiently long finite sequence of an everywhere dense sequence in $E$). 

We can then implement the Tchakaloff-Davis-Wilhelmsen result on a spline curved element $E$ by moment-matching in the following way:

\begin{itemize}
\item take the bivariate Chebyshev polynomial basis of a cartesian bounding box $B\supseteq E$ with center $(c_1,c_2)$ and side half-lengths $\ell_1$ and $\ell_2$
\begin{equation} \label{bivcheb}
f_{\alpha}(x,y)=T_{\alpha_1}\left(\frac{x-c_1}{\ell_1}\right)\,T_{\alpha_2}\left(\frac{y-c_2}{\ell_2}\right)\;,\;(x,y)\in B\;,\;\;\alpha=(\alpha_1,\alpha_2)\;,\;0\leq \alpha_1+\alpha_2\leq n\;,
\end{equation}
where $T_s(t)=\cos(s\arccos(t))$ are the Chebyshev polynomials of the first kind in $[-1,1]$ and $(x,y)\mapsto (\ell_1x+c_1,\ell_2y+c_2)$ is the standard affine transformation of the bounding box into the reference square $[-1,1]^2$; the $N=dim(\mathbb{P}_n)=(n+1)(n+2)/2$ couples $(\alpha_1,\alpha_2)$ are suitably (e.g. lexicographically) ordered

\item 
consider the $x$-primitives
$\mathcal{F}_\alpha(x,y)=\int{f_\alpha(x,y)\,dx}$ 
\begin{equation} \label{chebprim}
\mathcal{F}_\alpha(x,y)=\ell_1\,\mathcal{F}_{\alpha_1}\left(\frac{x-c_1}{\ell_1}\right)\,T_{\alpha_2}\left(\frac{y-c_2}{\ell_2}\right)\;,
\end{equation}
where (cf. \cite[\S 2.4.4]{MH02})
$$
\mathcal{F}_{\alpha_1}(t)=\frac{1}{2}\left(\frac{T_{\alpha_1+1}(t)}{\alpha_1+1}-\frac{T_{|\alpha_1-1|}(t)}{\alpha_1-1}\right)\;,\;\alpha_1\neq 1\;,\;\;\mathcal{F}_1(t)=\frac{T_2(t)}{4}
$$

\item compute the $N=(n+1)(n+2)/2$ Chebyshev moments as in (\ref{splinegauss})
$$
m_\alpha=\int_E{f_\alpha(x,y)\;dxdy}=\oint_{\partial E}{\mathcal{F}_\alpha(x,y)\,dy}
$$
\begin{equation} \label{chebmom}
=\sum_{i=1}^m\sum_{j=1}^{s_i}\sum_{h=1}^{\lceil dn/2\rceil}{\lambda _{i,j,h}\,\mathcal{F}_\alpha(S_{i}(\tau_{i,j,h}))\,y'_{i}(\tau_{i,j,h})}\;,\;\;\;0\leq \alpha_1+\alpha_2\leq n
\end{equation}

\item take a subsequence $\mathcal{H}_{N_k}$ of the Halton sequence of the bounding box $B\supseteq E$, e.g. $N_k=(1+\theta)N_{k-1}$, $k\geq 1$, for a suitable $\theta\in (0,1]$, and compute the sequence 
\begin{equation} \label{Hk}
H_{k}=\mathcal{H}_{N_k}\cap E=\{Z_1,\dots,Z_{M_k}\}
\end{equation} 
by the {\em in-domain algorithm} based on the spline boundary {\em crossing number} as implemented in \cite{SV21}; $N_0$ has to be chosen sufficiently large to ensure that $M_0\approx N_0\times area(E)/area(B)\geq N$, say $N_0>(1+\theta)N\times area(B)/area(E)$, where the element area can be computed 
by (\ref{splinegauss}) with $\mathcal{F}(x,y)=x$

\item compute the transposed Chebyshev-Vandermonde matrices $V_k^t=(f_\alpha(Z_\ell))_{\alpha,\ell}\in \mathbb{R}^{N\times M_k}$
and solve the NonNegative Least Squares problem (all vectors here are column vectors)
\begin{equation} \label{nnls}
\|V_k^tu_k-\{m_\alpha\}\|_2^2=\min_{0\leq u\in \mathbb{R}^{M_k}}\|V_k^tu-\{m_\alpha\}\|_2^2\;,\;\;k=0,1,2,\dots
\end{equation}
by the Lawson-Hanson active set method \cite{LH95}, until the residual goes below a given tolerance for a suitable $k^\ast$
(see \cite{SV15} for a discussion on the effect of a nonzero residual on the cubature error)

\item select the $M=N$ positive weights $\{w_s\}=\{u_{k^\ast}(\ell)>0\}$ and nodes $\{P_s\}=\{Z_\ell:u_{k^\ast}(\ell)>0\}$ of the final cubature formula.
\end{itemize}

This algorithm has some appealing features: {\em no subtessellation} of the element is needed, the final weights are {\em positive}, the nodal set is {\em inside the element} and has {\em low cardinality}. The main drawback is an {\em high computational cost}. The bulk is given by Lawson-Hanson algorithm applied at each step $k$, that by construction is able to select a sparse solution with  $M\leq N$ nonzero weights, but is based on a sequence of QR factorizations with column pivoting \cite{BG65} (in practice, a sequence of sequences of QR factorizations is performed). 
Even when an accelerated version of Lawson-Hanson is used, adopting the concept of {\em deviation maximization} instead of column pivoting \cite{DODM23,DM22,DMV20}, the method appears costly for an elementwise application within FEM or VEM on a single problem. The situation does not improve if moment-matching is implemented via linear programming instead of quadratic programming, as tested in \cite{PSV17}. On the other hand, the method becomes interesting when 
a given curved discretization mesh for a PDE has to be used several times, for example with {\em parametric or time-dependent problems}, since the cubature rule can be computed only once for each spline curved element.

\subsection{{\bf Spline-Chebyshev}}
The Spline-Chebyshev approach aims at reducing the computational cost of Spline-Tchakaloff, still keeping a low cardinality of the cubature rule. 
The method was inspired by the ``frugal'' approach to algebraic cubature on polytopal elements, recently proposed in \cite{LvPBD24} to accelerate the computation of stiffness and mass matrices. It is still a moment-matching and subtessellation-free method, where in order to get at the same time a low computational cost and a low cardinality of the rule, one accepts a priori to support the rule also at external points, namely at suitable points of a bounding box for the element, and to have some possible negative weights. 

In \cite{LvPBD24} this was done by polynomial interpolation at approximate Fekete points \cite{SV09-2} of the bounding box, that can be computed by a single QR factorization with column pivoting of a transposed Vandermonde matrix at a suitable grid of the box. The resulting weights are not all positive, in general, but empirical stability was observed. While the method had some manifest advantages (no triangulation of the polytopal element needed, QR factorization made only once on a reference box), it also presented the drawback of solving a linear system for each element, possibly ill-conditioned due to the adoption of the monomial basis to simplify moment computation by the recursive algorithm in \cite{AHP18}. 

In \cite{RSV26} we proposed a similar approach, based on {\em Chebyshev-like hyperinterpolation} instead of interpolation. We recall that Sloan hyperinterpolation \cite{Slo95} is an orthogonal projection on $\mathbb{P}_n$ with respect to an absolutely continuous measure $d\mu=w(P)\,dP$, where the scalar products in $L^2_{d\mu}$ are discretized by a cubature rule with exactness degree $2n$, that is a discretized truncated expansion by an orthonormal polynomial basis. 
Such a basis is here the {\em orthonormal Chebyshev basis} of 
a cartesian bounding box $B\supseteq E$ with center $(c_1,c_2)$ and side half-lengths $\ell_1$ and $\ell_2$
$$
\{g_\alpha(x,y)\}=\left\{\frac{1}{\sqrt{\ell_1\ell_2}}\,\hat{f}_\alpha(x,y)\right\}\;,\;\;\alpha=(\alpha_1,\alpha_2)\;,\;0\leq \alpha_1+\alpha_2\leq n\;,
$$
\begin{equation} 
\hat{f}_\alpha(x,y)=c_{\alpha_1}c_{\alpha_2}\,T_{\alpha_1}\left(\frac{x-c_1}{\ell_1}\right)\,T_{\alpha_2}\left(\frac{y-c_2}{\ell_2}\right)\;,\;\;\;(x,y)\in B\;,
\end{equation}
where $c_s$ is the  normalization constant for the orthogonal Chebyshev basis $\{T_s\}$ in $[-1,1]$. On the other hand, the cubature rule is given by the positive weights $\{u_s\}$ (scaled by $\ell_1\ell_2$) and affinely mapped nodes $\{(\xi_s,\eta_s)\}$ of the {\em near-minimal Morrow-Patterson-Xu rule} for the Chebyshev measure of the first kind in the reference square $[-1,1]^2$, with {\em polynomial exactness degree} $2n$ and cardinality 
\begin{equation} \label{MPX}
M=(n+1)(n+3)/2\;,\;n\;odd\;,\;\;\;\;M=(n+2)^2/2\;,\;n\;even\;,
\end{equation}
cf. \cite{CDMMV06,X25} and \cite{S25}. 

The method can be then implemented as follows:

\begin{itemize}
\item compute the Morrow-Patterson-Xu weights $\{u_s\}$ and nodes $\{(\xi_s,\eta_s)\}\subset [-1,1]^2$, $1\leq s\leq M$ (cf. \cite{S25})

\item compute the Chebyshev-Vandermonde matrix $\hat{V}=(c_{\alpha_1}c_{\alpha_2}\,T_{\alpha_1}(\xi_s)\,T_{\alpha_2}(\eta_s))_{s,\alpha}\in \mathbb{R}^{M\times N}$
\end{itemize}

then, {\em for every element} $E$

\begin{itemize}

\item make the first three steps (\ref{bivcheb})-(\ref{chebmom}) of Spline-Tchakaloff to compute the moments $\{\hat{m}_\alpha\}=\{c_{\alpha_1}c_{\alpha_2}\,m_\alpha\}$;

\item compute the $M$ nodes $\{P_s\}=\{(\ell_1\xi_s+c_1,\ell_2\eta_s+c_2)\}\subset B$ and weights $\{w_s\}=diag(u_s)\hat{V}\{\hat{m}_\alpha\}$ of the final cubature rule.
\end{itemize}

The proof that the rule is exact in $\mathbb{P}_n$ is an immediate consequence of the fact that $\hat{V}=(\hat{f}_\alpha(P_s))_{s,\alpha}$ and that the matrix $diag(\sqrt{u_s})\hat{V}$ is {\em orthogonal}, 
and thus 
\begin{equation} \label{chebw}
(\hat{f}_\alpha(P_s))^t\{w_s\}=\hat{V}^t\{w_s\}=\hat{V}^tdiag(u_s)\hat{V}\{\hat{m}_\alpha\}=\{\hat{m}_\alpha\}\;,
\end{equation}
i.e. the rule is exact on the polynomial basis $\{\hat{f}_\alpha(x,y)\}$. 

Notice that the first two steps are {\em element independent} and can be done once and for all. The bulk of the algorithm is now in the moment computation, because then the weights are obtained by {\em a single matrix-vector product} (with fixed matrix of relatively small size, 
since $M\approx N=(n+1)(n+2)/2$), followed  by a componentwise scaling 
of the resulting vector. Moreover, since {\em no matrix factorization} and {\em no system solving} is performed, we completely {\em avoid conditioning issues}. In addition, as proved in \cite{SV25} not only the stability parameter $\sigma$ 
in (\ref{sigma}) is bounded but also $\sigma\to 1$ as $n\to \infty$, 
so that the cubature rule is {\em almost optimally stable} even with some negative weights.

Such a combination of {\em no conditioning problem, low construction cost, stability and low cardinality}, make Spline-Chebyshev a valid competitor with the method proposed in \cite{LvPBD24} and a good choice in FEM and VEM computations, with a main drawback: like \cite{LvPBD24} the method relies on the possibility of sampling in a bounding box, i.e. also outside the element. This feature, shared by Spline-Gauss in nonconvex instances, is not a problem for PDEs with smooth coefficients. On the other hand, it becomes a limitation in cases where one may have coefficients which are jumping
between elements, or are not even defined outside the computational domain where
the PDE is to be numerically approximated. This is the reason 
that has led to the fourth method, Spline-Fekete, described in the next subsection.

\subsection{{\bf Spline-Fekete}}
The idea of constructing cubature rules by polynomial interpolation at Fekete-like points of complex-shaped integration domains, is not new (cf. e.g. \cite{GSV11}), but is here codified for the first time in the application to curved elements with spline boundary. As observed above, the main motivation is to provide implementers of curved polytopal FEM/VEM with algebraic cubature rules that have {\em low-cardinality internal nodal sets} and {\em acceptable stability features}, at a {\em low/moderate computational cost}. 

We recall that Fekete-like sets (also called discrete extremal sets) are finite point distributions on a compact set that {\em (nearly) maximize} the absolute value of the {\em Vandermonde determinant}, giving good stability and convergence properties to polynomial interpolation 
(Lebesgue constants with $\mathcal{O}(N)$ growth, 
but often much lower in practice) and to the correspondent interpolatory cubature. Algorithms for their extraction from polynomial meshes, that are 
sequences of low-cardinality finite norming sets for the uniform norm on polynomial spaces, have been deeply studied over the last two decades. One of the most relevant properties is that the discrete uniform probability measure supported at such discrete extremal sets converges (weakly) to the pluripotential theoretic equilibrium measure of the compact set; cf. e.g. \cite{BCLSV11,BDMSV10,SV09-2}. 

We focus here on {\em approximate Fekete points} of an element $E$, 
which can be computed by QR factorization with column pivoting of a transposed Vandermonde matrix at a suitable discretization of $E$. The Spline-Fekete algorithm can be summarized as follows: 

\begin{itemize}

\item compute the moments $\{m_\alpha\}$ and the Halton sets $H_k=\{Z_\ell\}_{1\leq\ell\leq M_k}\subset E$ as in Spline-Tchakaloff;
    
\item compute the transposed Chebyshev-Vandermonde matrices $V_k^t=(f_\alpha(Z_\ell))_{\alpha,\ell}\in \mathbb{R}^{N\times M_k}$
and solve the underdetermined moment-matching system
\begin{equation} \label{matching}
V_k^tu_k=\{m_\alpha\}\;,\;\;k=0,1,2,\dots
\end{equation}
by QR factorization with column pivoting \cite{BG65} (i.e. directly by the backslash operator in Matlab), until the stability parameter $\sigma=\sum_\ell{|u_k(\ell)|}/area(E)$ goes below a given threshold for a suitable $k^\ast$;

\item select the $M=N$ nonzero weights $\{w_s\}=\{u_{k^\ast}(\ell)\neq 0\}$ and nodes $\{P_s\}=\{Z_\ell:u_{k^\ast}(\ell)\neq 0\}$ of the final cubature formula.

\end{itemize}

We stress that the resulting nonzero weights are not all positive in general, but in all our numerical tests the negative ones turn out to be few and of small size. Though we have not a theoretical proof of this fact, but only the qualitative explanation that Fekete-like points are near-optimally distributed for polynomial interpolation, the final cubature rule has good stability properties for the degree range considered, 
with a stability parameter $\sigma$ staying not far from the optimal value 1.

\section{Numerical examples}
In this section we construct and test, for the purpose of illustration, the four cubature rules described above on two elements with spline boundary. The first is convex and has seven sides (cf. Figure 1), six linear and one curved, while the second is nonconvex (cf. Figure 4) and has five sides, four linear and one curved with a complicated shape. Both can be thought as the intersection of convex polygonal elements with a curved boundary. The curved sides are tracked by cubic splines 
with ten control points on the convex arc, seven control points on the concave arc, and not-a-knot additional constraints. 
All the tests have been performed by Matlab v2024b on an Intel Core Ultra 5 125H processor. 

The distribution of cubature nodes at exactness degree $n=10$ is displayed in Figures 1 and 4. 
Notice that by construction Spline-Tchakaloff, Spline-Fekete, and Spline-Gauss in the convex case, have internal nodes, whereas Spline-Chebyshev has by construction also external nodes and Spline-Gauss can have external nodes in nonconvex instances. The distribution of the corresponding weights are displayed in Figures 2 and 5, respectively. We recall that by construction Spline-Tchakaloff,
and Spline-Gauss in the convex case, have positive weights, whereas negative weights (but few and of small size) occur with Spline-Chebyshev, Spline-Fekete (where they turn out to be located 
at the boundary), and can occur for Spline-Gauss in nonconvex instances (unless it is possible to choose a suitable reference line as described in \cite{SV09}).  

The essential parameters of the four cubature rules are reported in Tables 1-2 and Figures 3 and 6. The figures show the distribution and the geometric mean of the relative integration errors for 100 trials a family of random polynomials with degrees $n=2,4,6,\dots,16$. Such geometric means deviate of at most two orders of magnitude with respect to machine precision. 

In Tables 1-2, at the same sequence of exactness degrees, we display the rule cardinality, the stability parameter $\sigma=\sum_s|w_s|/area(E)$, and the CPU time. By construction, Spline-Tchakaloff and Spline-Fekete have a cardinality equal to $N=dim(\mathbb{P}_n)=(n+1)(n+2)/2$. On the other hand the cardinality of Spline-Chebyshev is only a little higher than $N$, being that of the Morrow-Patterson-Xu rule of a bounding box as in (\ref{MPX}), whereas the cardinality of Spline-Gauss is much bigger, cf. (\ref{cardsplineG}). The computing time is on the order of $10^{-4}$ seconds for Spline-Gauss and Spline-Chebyshev at all degrees, whereas for Spline-Tchakaloff and Spline-Fekete varies from an order of $10^{-3}$ to an order of $10^{-2}$ seconds (Spline-Fekete being in any case at least 3-4 times faster at the highest degrees).

We can conclude by the following observations on possible selection of the rule 
among the four presented (see also the feature sketch in Table 3): 
\begin{itemize}

\item all the rules, even those with some negative weights, show excellent stability properties in the degree range considered (suitable for FEM and VEM), with a stability parameter equal or very close to 1;

\item if computing time is the key parameter, as in large scale 
applications, Spline-Gauss and Spline-Chebyshev are suitable choices, 
with a drawback of the former concerning large cardinality of the support (that could affect the construction time of stiffness and mass matrices), especially when the same finite element discretization has to 
be used for different equations on a given domain (e.g. for parametric or time-dependent problems); 

\item if only internal nodes are desired, irrespectively of the element shape, then one should choose 
either Spline-Tchakaloff or Spline-Fekete (faster), for example in cases where one may have coefficients which are jumping
between elements, or are not even defined outside of the computational domain where
the PDE is to be numerically approximated;

\item if theoretically guaranteed stability at any degree is the desired feature, all the methods but Spline-Fekete 
can be used, with Spline-Tchakaloff providing also positive weights, as well as Spline-Gauss in convex instances.

\end{itemize}

\begin{figure}[h]
	\centering
	\hspace{-1cm}
    
	\includegraphics[scale=0.0595,clip,valign=t]{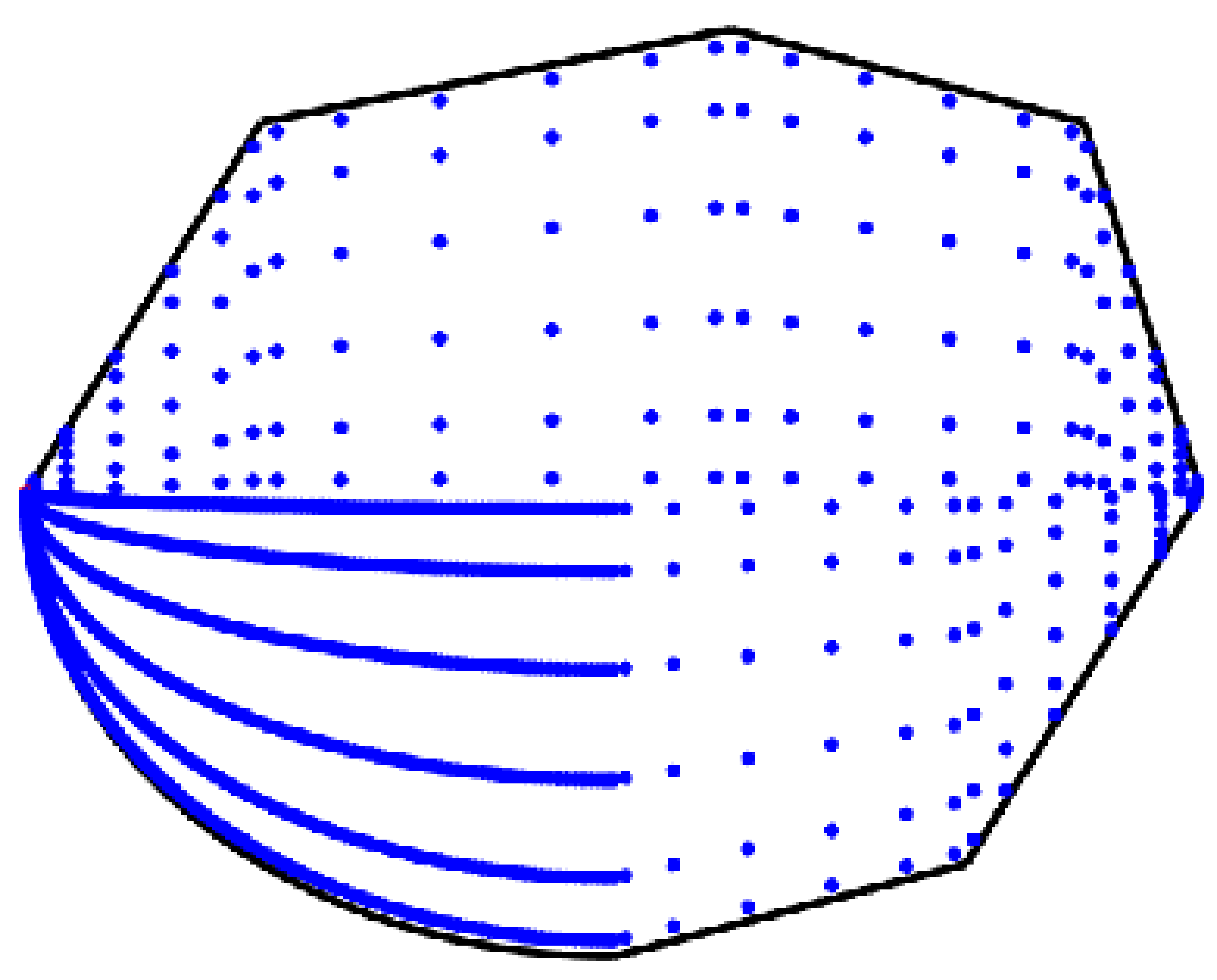}
	\hspace{1.5cm}
	\includegraphics[scale=0.50,clip,valign=t]{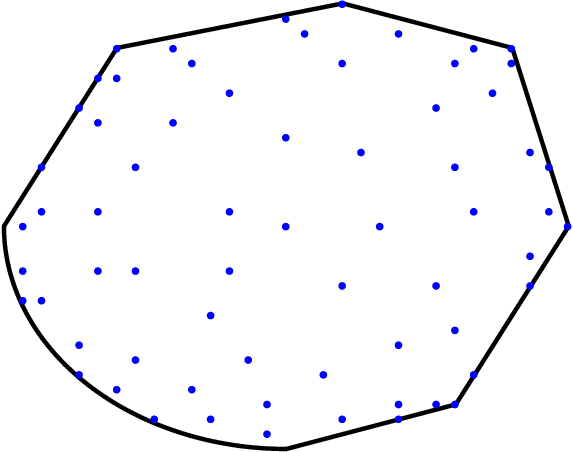}
    \vskip1cm
    \includegraphics[scale=0.50,clip,valign=t]{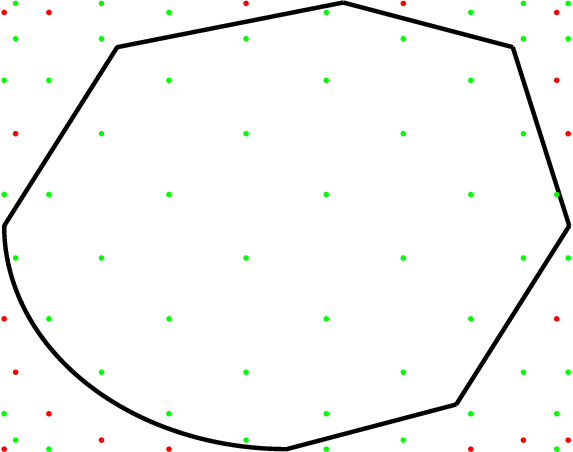}
	\hspace{1.5cm}
    \includegraphics[scale=0.50,clip,valign=t]{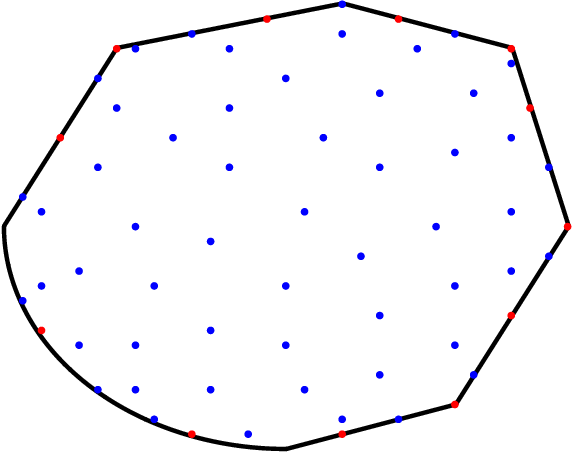}
	\caption{Nodes of the four rules with polynomial degree of exactness $n=10$ on a convex spline element $E$: Spline-Gauss (top-left, 684 nodes), Spline-Tchakaloff (top-right, 66 nodes), Spline-Chebyshev (bottom-left, 72 nodes), Spline-Fekete (bottom-right, 66 nodes). Blue dots: nodes with positive weights; red dots: nodes with negative weights.}
	\label{fig1:points}
\end{figure}

\begin{figure}[h]
	\centering
	\hspace{-0.35cm}
	\includegraphics[scale=0.2805,clip,valign=t]{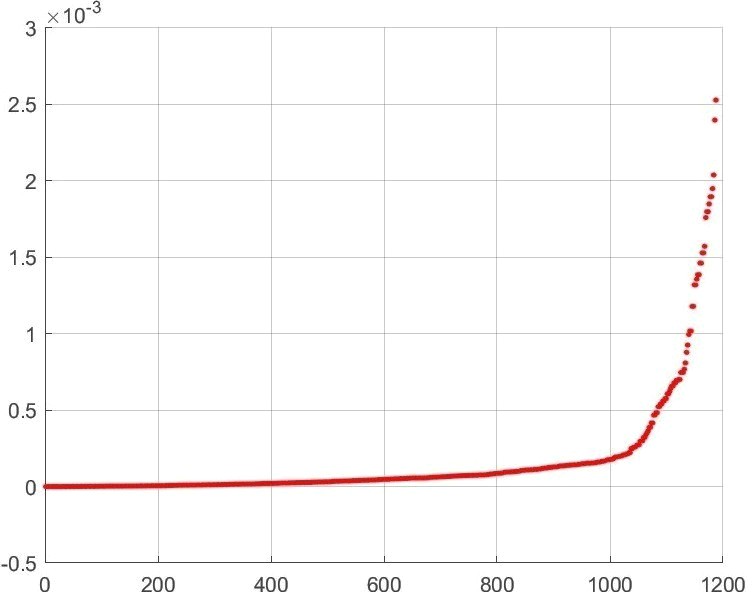}
	\hspace{0.275cm}
	\includegraphics[scale=0.45,clip,valign=t]{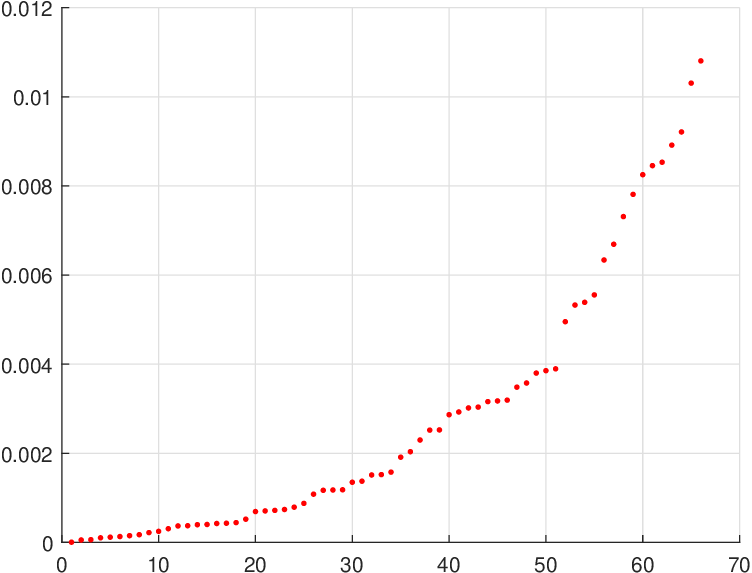}
    \vskip0.5 cm
    \includegraphics[scale=0.45,clip,valign=t]{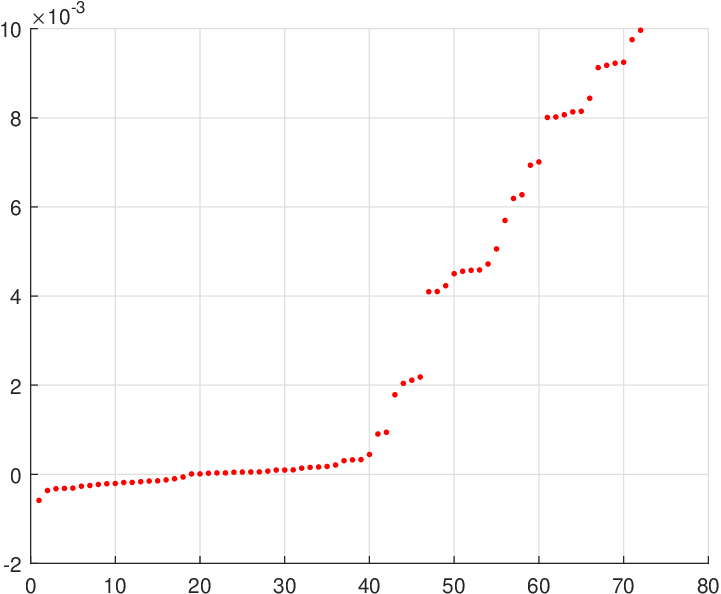}
	\hspace{0.6cm}
	\includegraphics[scale=0.45,clip,valign=t]{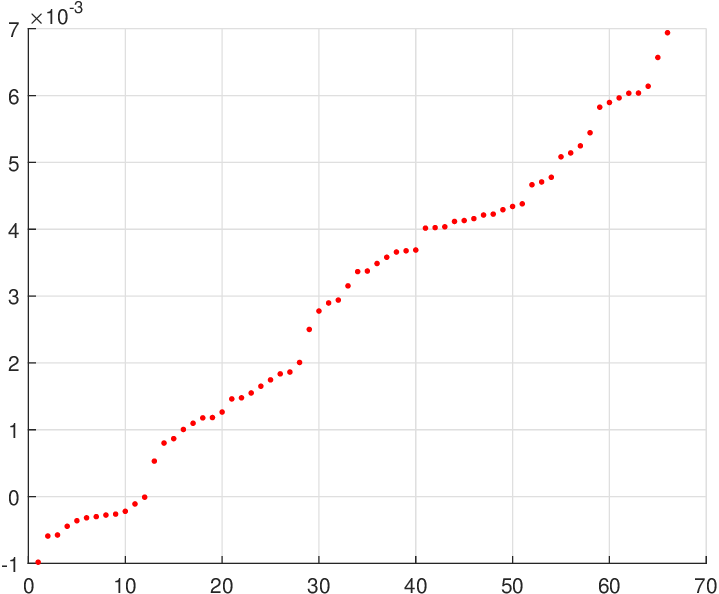}
	\caption{Weights magnitude distribution of the four rules with polynomial degree of exactness $n=10$ on the convex spline element $E$ of Figure \ref{fig1:points}: Spline-Gauss (top-left), Spline-Tchakaloff (top-right), Spline-Chebyshev (bottom-left), Spline-Fekete (bottom-right).}
	\label{fig1:weights}
\end{figure}

\begin{figure}[h]
	\centering
	\hspace{-1cm}
    
	\includegraphics[scale=0.45,clip,valign=t]{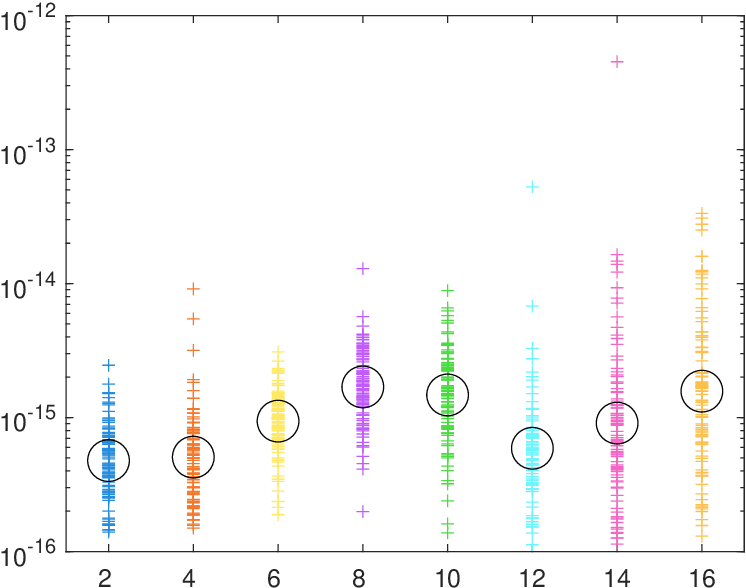}
	\hspace{0.6cm}
	\includegraphics[scale=0.45,clip,valign=t]{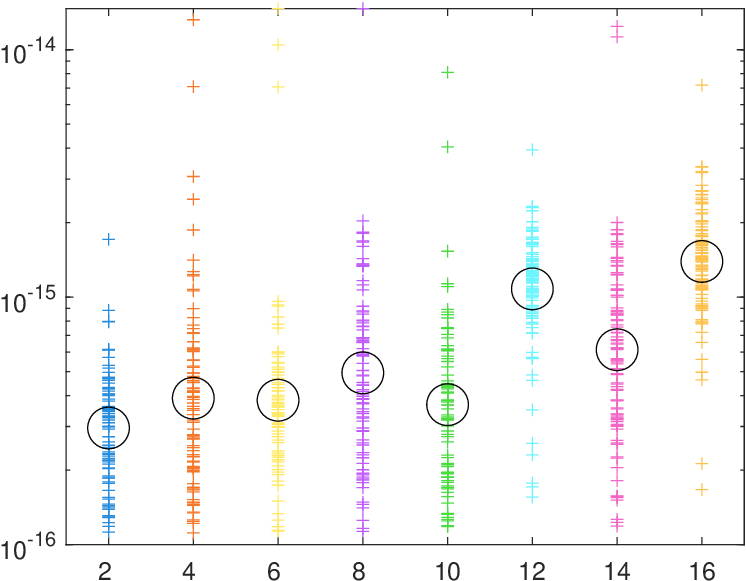}
    \vskip0.5 cm
    \includegraphics[scale=0.45,clip,valign=t]{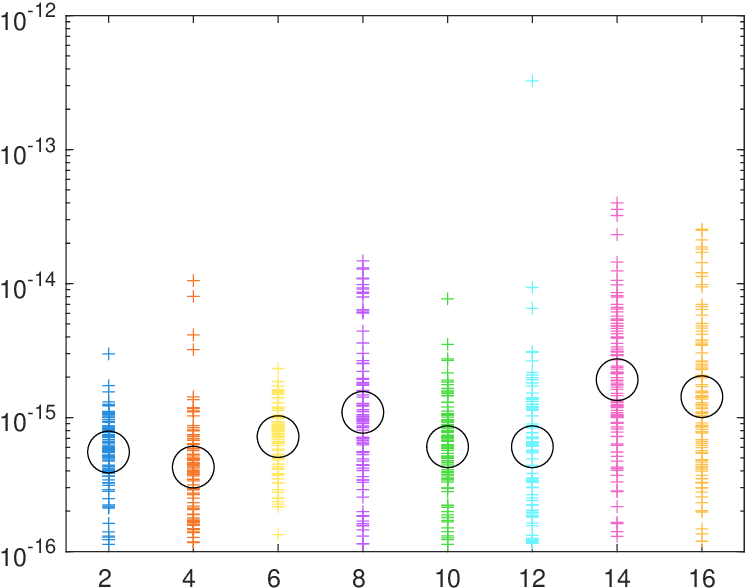}
	\hspace{0.6cm}
	\includegraphics[scale=0.45,clip,valign=t]{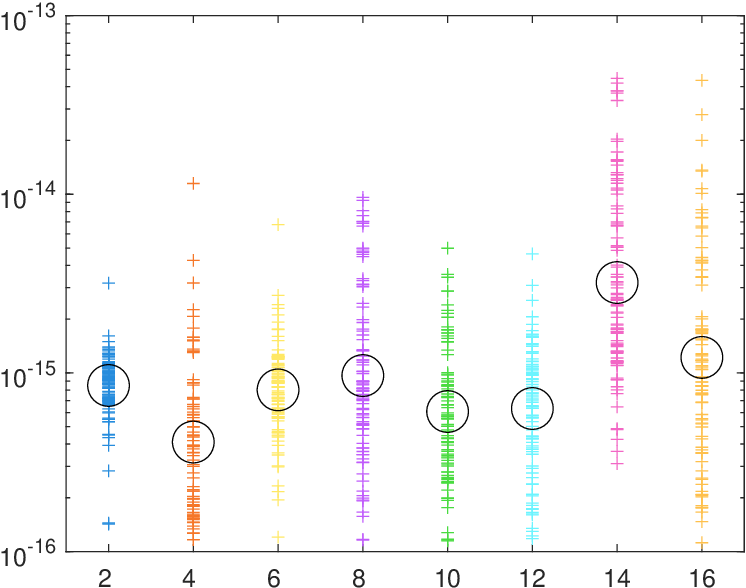}
	\caption{Geometric means of the integration relative  errors for 100 trials of random polynomials $(c_0+c_1 x+c_2 y)^n$ on the convex spline element $E$ of Figure \ref{fig1:points}: Spline-Gauss (top-left), Spline-Tchakaloff (top-right), Spline-Chebyshev (bottom-left), Spline-Fekete (bottom-right).}
	\label{fig1:geomeans}
\end{figure}

\begin{table}
	\centering
	\begin{tabular}{|c||c|c|c||c|c|c||c|c|c||c|c|c|}  
      \hline  
      & \multicolumn{3}{c||}{Spline-Gauss} & \multicolumn{3}{c||}{Spline-Tchakaloff} & \multicolumn{3}{c||}{Spline-Chebyshev} & \multicolumn{3}{c|}{Spline-Fekete}\\
      \hline  
       deg $n$& card & stab & cpu & card & stab & cpu & card & stab & cpu & card & stab & cpu\\
      \hline  
      2 & 132 & 1.00 & 7.5e-4 & 6 & 1.00 & 2.0e-3 & 8 & 1.21 & 5.5e-4 & 6 & 1.03 & 2.0e-3 \\
      4 & 297 & 1.00 & 7.1e-4& 15 & 1.00 & 2.2e-3& 18 & 1.08 & 6.8e-4& 15 & 1.10 & 2.1e-3\\
      6 & 528 & 1.00 & 8.4e-4 & 28 & 1.00 & 3.9e-3 & 32 & 1.06 & 6.3e-4 & 28 & 1.00 & 3.5e-3 \\
      8 & 825 &  1.00 & 7.2e-4 & 45 & 1.00 & 6.1e-3 & 50 & 1.06 & 6.3e-4 & 45 & 1.05 & 2.9e-3 \\
      10 & 1188 & 1.00 & 7.6e-4 & 66 & 1.00 & 1.3e-2 & 72  & 1.04 & 6.9e-4 & 66 & 1.05 & 4.6e-3 \\
      12 & 1617 &  1.00 & 7.8e-4 & 91 & 1.00 & 3.2e-2 & 98 & 1.05 & 7.5e-4 & 91 & 1.01 & 8.0e-3 \\
      14 & 2112 & 1.00 & 8.3e-4 & 120 & 1.00 & 5.1e-2 & 128 & 1.04 & 8.2e-4 & 120 & 1.02 & 1.5e-2\\
      16 & 2673 & 1.00 & 1.0e-3 & 153 & 1.00 & 6.8e-2 & 162 & 1.04 & 8.5e-4 & 153 & 1.06 & 2.6e-2\\
\hline
\end{tabular}
\caption{Essential parameters of the four rules on the convex spline element $E$ of Figure \ref{fig1:points}; rule cardinality, the stability parameter $\sigma=\sum_s|w_s|/area(E)$, cputime in seconds.}
\label{tab:convex}
\end{table}

\begin{figure}[h]
	\centering
	\hspace{-1cm}
    
	\includegraphics[scale=0.50,clip,valign=t]{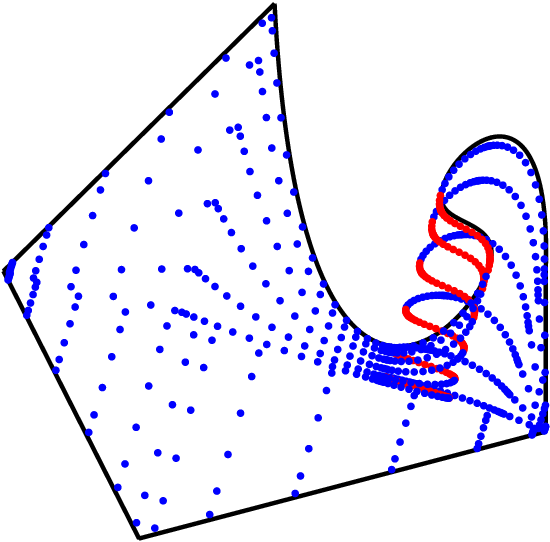}
	\hspace{1.5cm}
	\includegraphics[scale=0.50,clip,valign=t]{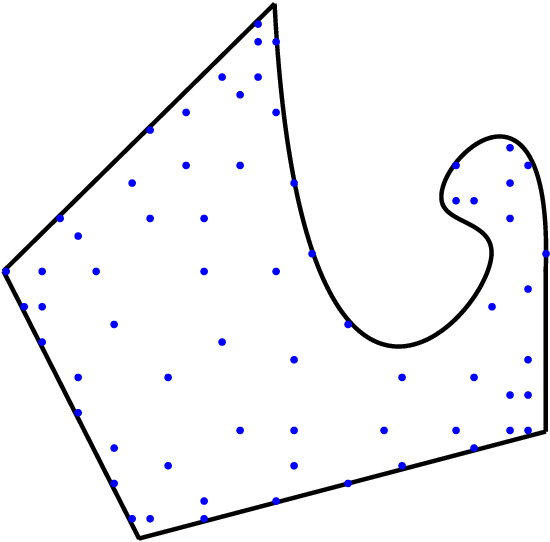}
    \vskip1cm
    \includegraphics[scale=0.50,clip,valign=t]{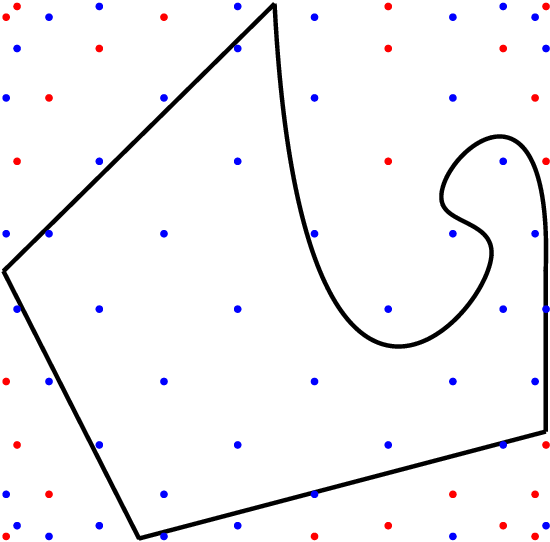}
	\hspace{1.5cm}
	\includegraphics[scale=0.50,clip,valign=t]{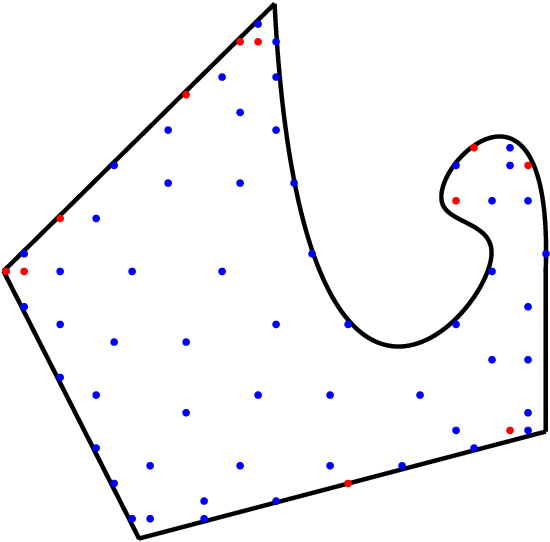}
	\caption{Nodes of the four rules with polynomial degree of exactness $n=10$ on a concave spline element $E$: Spline-Gauss (top-left, 684 nodes), Spline-Tchakaloff (top-right, 66 nodes), Spline-Chebyshev (bottom-left, 72 nodes), Spline-Fekete (bottom-right, 66 nodes). Blue dots: nodes with positive weights; red dots: nodes with negative weights.}
	\label{fig2:points}
\end{figure}

\begin{figure}[h]
	\centering
	\hspace{-1cm}
    
	\includegraphics[scale=0.45,clip,valign=t]{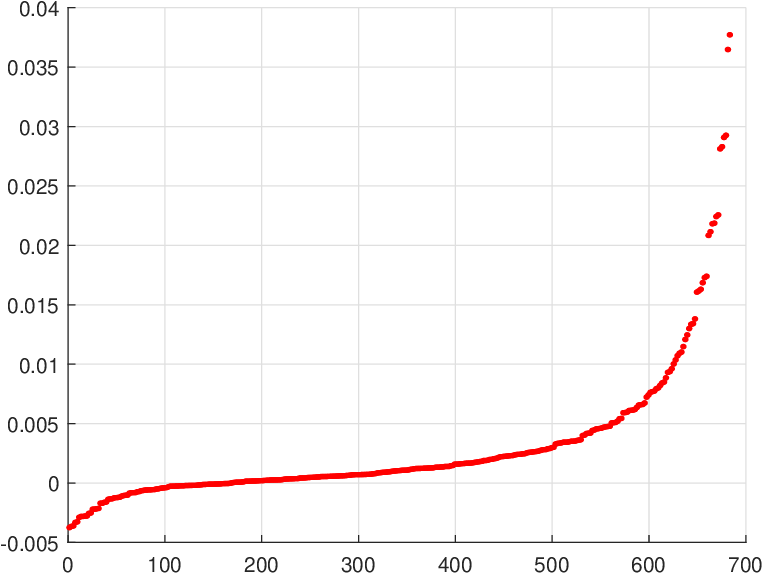}
	\hspace{0.6cm}
	\includegraphics[scale=0.45,clip,valign=t]{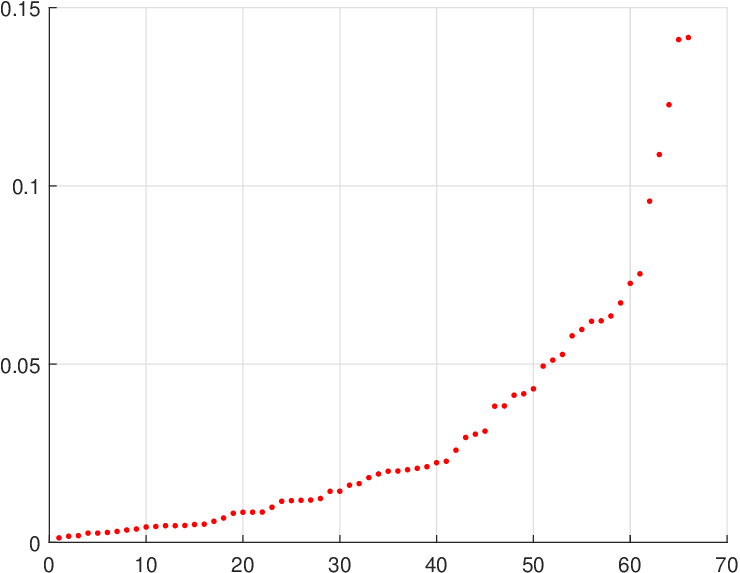}
    \vskip0.5 cm
    \includegraphics[scale=0.45,clip,valign=t]{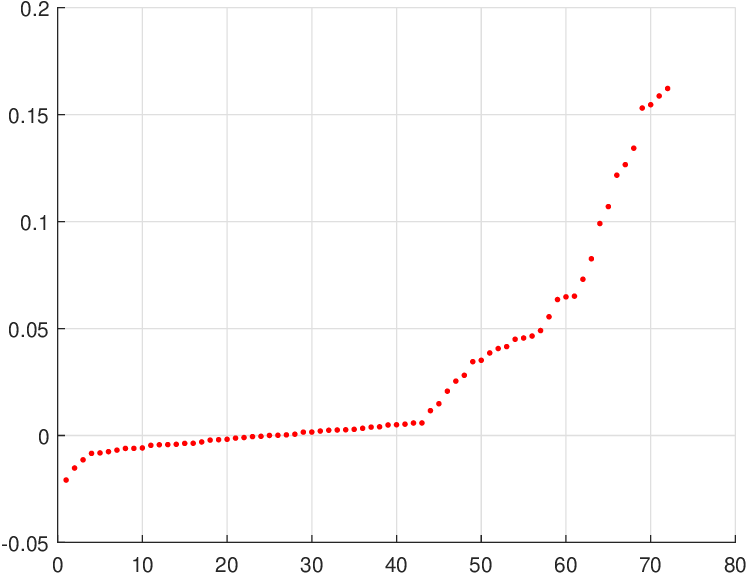}
	\hspace{0.6cm}
	\includegraphics[scale=0.45,clip,valign=t]{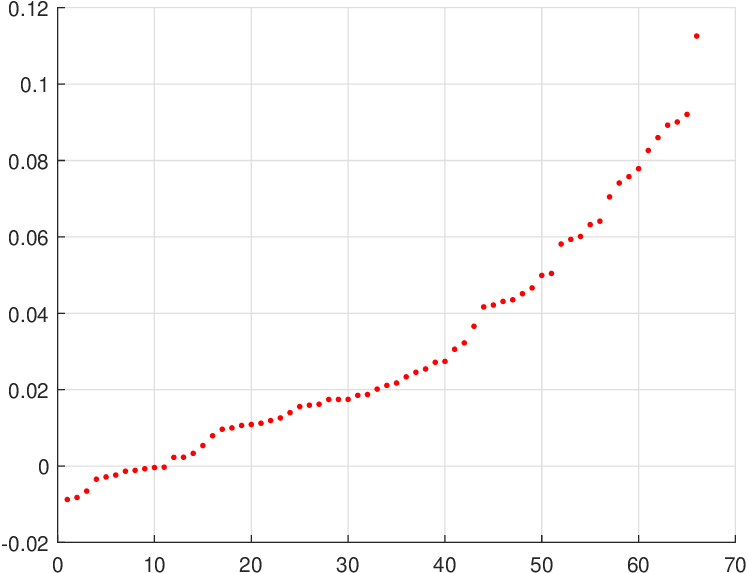}
	\caption{Weights magnitude distribution of the four rules with polynomial degree of exactness $n=10$ on the concave spline element $E$ of Figure \ref{fig2:points}: Spline-Gauss (top-left), Spline-Tchakaloff (top-right), Spline-Chebyshev (bottom-left), Spline-Fekete (bottom-right).}
	\label{fig2:weights}
\end{figure}

\begin{figure}[h]
	\centering
	\hspace{-1cm}
    
	\includegraphics[scale=0.45,clip,valign=t]{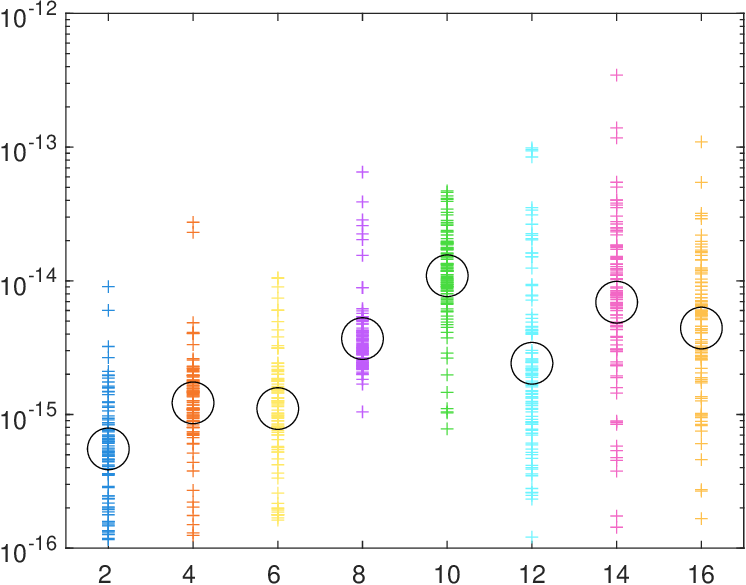}
	\hspace{0.6cm}
	\includegraphics[scale=0.45,clip,valign=t]{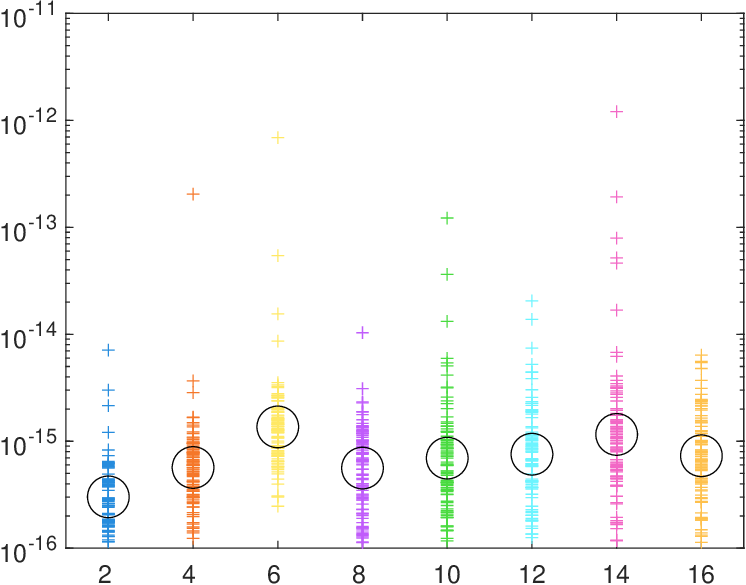}
    \vskip0.5 cm 
    \includegraphics[scale=0.45,clip,valign=t]{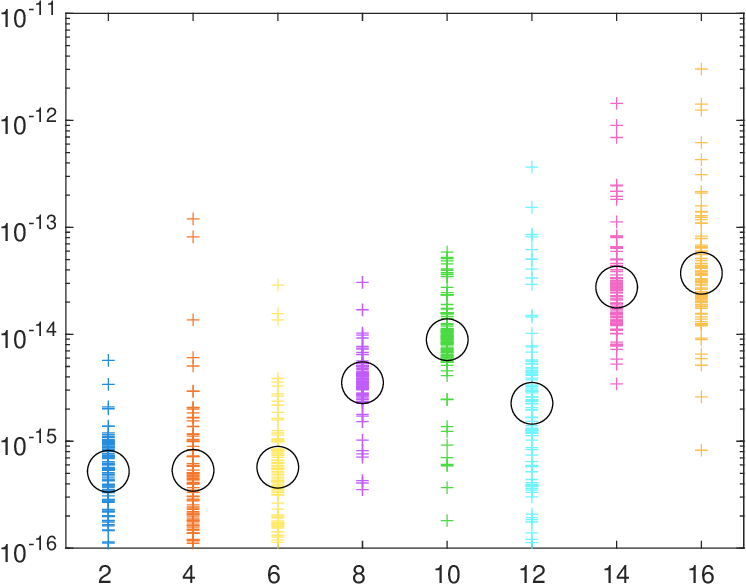}
	\hspace{0.6cm}
	\includegraphics[scale=0.45,clip,valign=t]{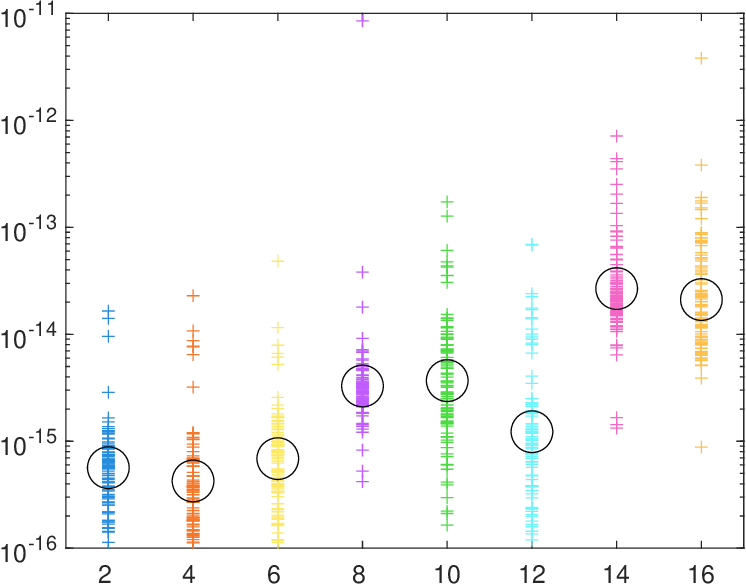}
	\caption{Geometric means of the integration relative  errors for 100 trials of random polynomials $(c_0+c_1 x+c_2 y)^n$ on the concave spline element $E$ of Figure \ref{fig2:points}: Spline-Gauss (top-left), Spline-Tchakaloff (top-right), Spline-Chebyshev (bottom-left), Spline-Fekete (bottom-right).}
	\label{fig2:geomeans}
\end{figure}

\begin{table}
	\centering
	\begin{tabular}{|c||c|c|c||c|c|c||c|c|c||c|c|c|}  
      \hline  
      & \multicolumn{3}{c||}{Spline-Gauss} & \multicolumn{3}{c||}{Spline-Tchakaloff} & \multicolumn{3}{c||}{Spline-Chebyshev} & \multicolumn{3}{c|}{Spline-Fekete}\\
      \hline  
       deg $n$& card & stab & cpu & card & stab & cpu & card & stab & cpu & card & stab & cpu\\
      \hline  
      2 & 76 & 1.16 & 7.9e-4 & 6 & 1.00 & 1.5e-3 & 8 & 1.24 & 4.2e-4 & 6 & 1.00 & 1.4e-3 \\
      4 & 171 & 1.16 & 7.8e-4& 15 & 1.00 & 1.7e-3& 18 & 1.20 & 4.9e-4& 15 & 1.16 & 1.5e-3\\
      6 & 304 &  1.16 & 7.7e-4 & 28 & 1.00 & 2.6e-3 & 32 & 1.15 & 4.4e-4 & 28 & 1.05 & 1.8e-3 \\
      8 & 475 &  1.16 & 7.7e-4 & 45 & 1.00 & 1.5e-2 & 50 & 1.10 & 4.6e-4 & 45 & 1.17 & 3.1e-3 \\
      10 & 684 & 1.16 &  8.1e-4 & 66 & 1.00 & 1.3e-2 & 72 & 1.13 & 4.8e-4 & 66 & 1.04 & 3.4e-3 \\
      12 & 931 &  1.16 & 8.7e-4 & 91 & 1.00 & 4.0e-2 & 98 & 1.08 & 5.5e-4 & 91 & 1.10 & 7.2e-3 \\
      14 & 1216 & 1.16 & 8.7e-4 & 120 & 1.00 & 6.9e-2 & 128 & 1.08 & 6.7e-4 & 120 & 1.22 & 1.3e-2\\
      16 & 1539 & 1.16 & 9.7e-4 & 153 & 1.00 & 8.5e-2 & 162 &  1.09 & 6.5e-4 & 153 & 1.13 & 2.1e-2\\
\hline
\end{tabular}
\caption{Essential parameters of the four rules on the concave spline element $E$ of Figure \ref{fig2:points}; rule cardinality, the stability parameter $\sigma=\sum_s|w_s|/area(E)$, cputime in seconds.}
\label{tab:concave}
\end{table}

\begin{table}

	\centering
	\begin{tabular}{||c||c||c||c||c|}  
      \hline  
   main features & {Spline-Gauss} & {Spline-Tchakaloff} & {Spline-Chebyshev} & {Spline-Fekete}\\
      \hline  
      sub-triangulation free & yes & yes & yes & yes\\
    \hline
    matrix factoriz. free & yes & no & yes & no\\
    \hline
    computing time  & low & high & low & medium\\
    \hline
    cardinality  & high & low & low & low\\
    \hline
    only internal nodes & no (yes, $E$ convex) & yes & no & yes\\\hline
    weights positivity & no (yes, $E$ convex) & yes & no & no \\\hline
    theoretical stability & yes & yes & yes & no \\\hline
    empirical stability & yes & yes & yes & yes \\
     \hline
\end{tabular}
\caption{Feature sketch of the four rules on a general or convex element $E$.}
\label{tab:sketch}
\end{table}

\vskip0.5cm 
\noindent
{\bf Acknowledgements.} 

Work partially supported by the DOR funds of the University of Padova.  
This research has been accomplished within the Community of Practice 
``Green Computing" of the Arqus European University Alliance, towards sustainable computing for massive simulation and design by finite elements, 
and within the INdAM-GNCS and SIMAI Italian research communities.

\end{document}